\documentclass[12pt]{article}

\usepackage[numbers,sort&compress]{natbib} 
\usepackage{color}
\pdfoutput=1

\usepackage{amsfonts}
\usepackage{mathrsfs} 
\usepackage{amsfonts,amssymb}
\usepackage{exscale}

\usepackage{comment}

\usepackage{amsgen,amsmath,amstext,amsbsy,amsopn,amsfonts,amssymb,pifont, epsfig}
\usepackage{graphicx}

\newtheorem{theorem}{Theorem}[section]

\newtheorem{definition}[theorem]{Definition}
\newtheorem{cor}[theorem]{Corollary}

\newtheorem{example}{Example}
\newtheorem{remark}[theorem]{Remark}

\newcommand{\be}{\begin{equation}}
\newcommand{\ee}{\end{equation}}
\newcommand{\bea}{\begin{eqnarray}}
\newcommand{\eea}{\end{eqnarray}}
\newcommand{\beaa}{\begin{eqnarray*}}
\newcommand{\eeaa}{\end{eqnarray*}}
\newcommand{\bei}{\begin{itemize}}
\newcommand{\eei}{\end{itemize}}
\newcommand{\bee}{\begin{enumerate}}
\newcommand{\eee}{\end{enumerate}}
\newcommand{\bi}{\begin{itemize}}
\newcommand{\ei}{\end{itemize}}
\newcommand{\beq}{\begin{eqnarray*}}
\newcommand{\eeq}{\end{eqnarray*}}
\newcommand{\beqn}{\begin{eqnarray}}
\newcommand{\eeqn}{\end{eqnarray}}

\newcommand{\ignore}[1]{}{}

\def\FF{\mathcal{F}}

\def\qed{{\hfill $\Box$ \bigskip}}

\def\P{{\mathbb P}}

\def\R{{\mathbb R}}

\def\<{\langle}
\def\>{\rangle}

\def\TT{\mathcal{T}}

\def\R{\mathbb R}

\def\FF{\mathcal{F}}

\numberwithin{equation}{section}

\usepackage{color}
\definecolor{webgreen}{rgb}{0,.5,0}
\definecolor{webbrown}{rgb}{.8,0,0}
\definecolor{emphcolor}{rgb}{0.5,0.95,0.95}
\usepackage{hyperref}
\hypersetup{%
          colorlinks=true,
          linkcolor=webbrown,
          filecolor=webbrown,
          citecolor=webgreen,
          urlcolor=blue,
          breaklinks=true}

\begin{document}

\title{\bf Stochastic Periodic Solutions for Newtonian Systems via Lyapunov Method}

\author{Junxia Duan, Jifa Jiang$\thanks{The corresponding author. \newline \indent \quad Email: duanjunxia0802@163.com, jiangjf@htu.edu.cn, xujiescu@163.com}$, Jie Xu
\\ {\small School of Mathematics and Statistics, Henan Normal University, Xinxiang, Henan 453007, China}}

\maketitle

\begin{abstract}
This paper establishes an existence theory for  periodic probability solutions to both time periodically and stochastically forced Newton's equation, where the friction matrix is the Hessian of a twice continuously differentiable friction function. Employing the periodic probability solution theory established by Ji et al \cite{Qi1,Qi2} and the Lyapunov method, we prove its  existence under the assumptions that both the friction and potential functions tend to positive infinity at infinity, and their gradient inner product grows at least like an even-power polynomial. Provided the potential grows slightly faster at infinity, this result persists under bounded friction matrix perturbations. This largely confirms the stochastic Levinson conjecture proposed in \cite{Liyong}  with much less restrictions on frictions and potentials and resolves the open problem and the interesting problem posed in \cite [P. 342]{Liyong}.

Periodic probability solutions exist in two key cases:
\begin{enumerate}
  \item {\bf Polynomial Case}: Both functions are polynomials with positive definite leading terms (a minimal condition for dissipation).

  \item {\bf Plasma Physics Case}: The potential force is bounded, while the friction function has sufficiently high growth (as in many plasma physics models).
\end{enumerate}
\end{abstract}

\smallskip\noindent
{\bf Keywords}: Stochastic time-periodic Newtonian systems; Stochastic Levinson's conjecture; Damping Hamiltonian systems; Periodic probability solution; Lyapunov's method

\section{Introduction}\label{S:1}
In physics, it has been a fundamental principle that a Newtonian system
\begin{equation}\label{DNewtPF}
\ddot{x} + A(x,\dot{x},t)\dot{x} + \nabla V(x)=e(t)
\end{equation}
with friction must have a periodic solution. That is, the system (\ref{DNewtPF}) admits a periodic solution if the friction matrix $A(x,\dot{x},t)$ is symmetric, positive definite and $A(x,\dot{x},t), e(t)$ are periodic in time $t$.
Mathematically in 1940s, Levinson \cite{Levinson} conjectured  the principle is true, which brought important dissipativeness theory (see \cite{LiWY,Burton,Yoshizawa,Zabreiko} and references therein). For deterministic Newtonian systems, Li et al. \cite{LiWY} verified Levinson's conjecture by \cite [Theorem 4.2, P. 2619]{LiWY}:
\begin{theorem}\label{LiWY}
Let $V\in C^2(\mathbb{R}^n)$, $e(t)$  and $A(x,\dot{x},t)$ be continuous and periodic in time $t$ with period $\TT$ and $A(x,\dot{x},t)$ locally Lipschitz continuous with respect to $(x, \dot{x})$. Suppose that
\begin{description}
  \item{\rm (i)} $A(x,\dot{x},t)$ is symmetric;

  \item{\rm (ii)} there exist positive constants $\alpha,\eta, \gamma$ and $ l$ such that
 \begin{eqnarray}\label{e:H02}
  A(x,\dot{x},t)\ge \alpha I_n;\quad x^{\top}\nabla V(x) \ge \eta \|x\|^2 \, \vee \, \gamma \, x^{\top}A(x,\dot{x},t) \,x  \qquad\hbox{for}~\ \|x\|\ge l
  \end{eqnarray}
  where $I_n$ is the $n\times n$ identity matrix and $\top$ denotes the transposition.
\end{description}
Then system {\rm(\ref{DNewtPF})} admits at least a $\TT$-periodic solution.
\end{theorem}
 \eqref{e:H02} implies that for every fixed $x$, $A(x,\dot{x},t)$ is bounded with respect to $\dot{x}$. That is, the unboundedness of the friction stems from the variable $x$, which warrants greater attention to the situation of unbounded friction. Note that condition \eqref{e:H02} implies the potential's growth rate exceeds that of the friction by more than 2.  However, this does not hold for the well-known van der Pol equation under time-periodic forcing.

In many contexts, systems are subject to stochastic noise and time-periodic recurrence, which can be modeled as stochastic time-periodic Newtonian systems. From a statistical perspective, Jiang, Li, and Yang \cite{Liyong} utilized the concept of {\it periodicity in distribution} \textemdash~a natural framework for describing periodic phenomenon in stochastic processes. Building on this, they proposed the following stochastic analogue of Levinson's conjecture.
\medskip \\
{\it {\bf Conjecture:} Stochastic time-periodic Newtonian system
 \begin{eqnarray}\label{SNewtPF}
\ddot{x} + A(x,\dot{x},t)\dot{x} + \nabla V(x)=\dot{B}_t+e(t)
\end{eqnarray}
admits $\TT$-periodic solution in distribution
if this system is dissipative}.
\medskip \\
Using the Wong-Zakai approximation method, Lyapunov's method, and Horn's fixed point theorem, the authors proved this conjecture holds if the potential function $V(x)$ exhibits at most quadratic growth under condition (ii) (see \cite[Theorem 3.1]{Liyong}). However, they noted that their result cannot be applied to the equation
\begin{equation}\label{V4}
\ddot{x} + \dot{x} + x^3=\dot{B}_t+ \sin t
\end{equation}
(see \cite[Section 4, P. 342]{Liyong}).  Their numerical simulations indicate that Eq.~\eqref{V4} does admit a periodic solution in distribution, leaving this case as an open problem. The authors further suggested that stochastic Levinson's conjecture with more complex coefficients is an interesting problem for future research.

Since then, the stochastic Levinson conjecture has been widely studied. Jiang and Li \cite{LiJ} established sufficient conditions for the conjecture to hold, assuming the potential function exhibits at most quadratic growth and that the friction matrix is independent of $x, \dot{x}$. Ji et al. \cite{Qi1} studied the long-term behavior in the sense of distribution for stochastic differential equations (SDEs) with irregular coefficients using time-periodic Fokker-Planck equations. Instead of searching {\it periodic solutions in distribution} of SDEs, they worked with  the {\it periodic probability solutions} of the associated Fokker-Planck equation, and established the existence theory of a periodic probability solution, which is applicable to both SDEs with locally Lipschitz continuous coefficients and SDEs with Sobolev coefficients.   Subsequently in \cite{Qi2}, they demonstrated the convergence of these solutions for SDEs with less regular coefficients. As an application, both papers \cite{Qi1, Qi2} provided sufficient conditions for the stochastic Levinson conjecture to hold, leveraging the Lyapunov function constructed by Wu \cite{Wu}.
 For the theories of more complex time-recurrence solutions in distribution, we can refer to \cite{Liu1, Liu2, Liu3, Zhao3, Lu3} and their many references therein.

Notably, system (\ref{SNewtPF}) encompasses the time-periodic Langevin equation
\begin{equation}\label{PLangrvin}
\ddot{x} + \gamma(t)\, \dot{x} + \nabla_x V(x,t) = \dot{B}_t
\end{equation}
and time-periodic stochastic damping Hamiltonian systems (see \cite{Wu}). Feng, Zhao, and Zhong \cite[Theorem 4.17]{Zhao} further proved that the Langevin equation with constant positive friction admits a unique periodic solution in distribution.

 We point out that under local Lipschitz continuity of the SDEs' coefficients, the concept of a periodic solution in distribution (Jiang et al. \cite{Liyong, LiJ}) is subsumed by the concept of a periodic probability solution (Ji et al. \cite{Qi1,Qi2}). We refer the reader to Section 2 for details. This motivates us to investigate the stochastic Levinson conjecture in the broader sense by establishing the existence of periodic probability solutions  for stochastic time-periodic Newtonian systems with either relaxation oscillating friction, or uniformly positive friction, or beyond.

We first present a general existence criterion for periodic probability solutions to stochastic time-periodic Hessian-driven friction Newtonian systems  in which either the potential force or the friction exhibits at least polynomial growth at infinity:
\begin{equation}\label{e:2.11}
\begin{cases}
   dx_t= y_tdt, \\
    dy_t=-\big[D^2 F(x_t) \, y_t +\nabla_{x} V(x_t,t)+ E(x_t,y_t,t)\big] \,dt+\Sigma(x_t,y_t,t) \, d{B}_t,
\end{cases}
\end{equation}
which is stated as follows.

\begin{theorem}\label{UF}
The system \eqref{e:2.11} admits at least a $\TT$-periodic probability solution if it satisfies Hypotheses {\rm {\bf (H1)}-{\bf\bf(H5)}}. Here
\begin{description}
  \item{\rm \bf (H1)} All coefficients in system \eqref{e:2.11} are continuous in their respective domains and locally Lipschitz continuous in $(x,y)$. Let $F\in C^2(\mathbb{R}^n), V(x,t)\in C_\TT^{2,1}(\mathbb{R}^n  \times \R)$, and $E(x,y,t)$, $\Sigma(x,y,t)$  also be periodic in time $t$ with period $\TT$. Additionally, $V(x,t)$ is lower bounded.

  \item{\rm \bf (H2)} There exists a constant $a>0$ such that
  \begin{equation*}\label{Pinfinity}
  \lim_{\|x\|\rightarrow +\infty} \Big[V(x,t)+aF(x)-\frac{a^2}{2}\|x\|^2 \Big]=+\infty, \qquad \hbox{uniformly in}~ t\in\mathbb{R}.
  \end{equation*}

  \item{\rm \bf (H3)} There exist positive constants $b, M$ and $ m$ such that the friction-potential gradient inner product satisfies
  \begin{equation*}
  -\langle \nabla_x V(x,t), \nabla F(x)-ax\rangle \le M-b\|x\|^{2m}, \qquad~ \forall (x,t)\in \mathbb{R}^n\times \mathbb{R}.
  \end{equation*}

  \item{\rm \bf (H4)} $E(x,y,t)$ is bounded by $e$ and
      \begin{equation}\label{e}
      \lim_{\|x\|\rightarrow +\infty}\frac{e \|\nabla F(x)-ax\|}{\|x\|^{2m}}=0.
      \end{equation}

  \item{\rm \bf (H5)} There exist positive constants $c_1, c_2, M_1$ and $M_2$ with $(c_1+c_2)<1$ such that
  \begin{equation*}
  |V_t(x,t)|\le c_1 \, b\|x\|^{2m}+M_1, \qquad\quad\qquad\qquad\qquad~~ \forall (x,t)\in \mathbb{R}^n\times \mathbb{R},
  \end{equation*}
  \begin{equation*}
  \| \Sigma\|_2^2:={\rm Tr}[(\Sigma \Sigma^{\top}) (x,y,t)]\le 2c_2 \, (a\|y\|^2 +b\|x\|^{2m})+M_2, \quad \forall (x,y,t)\in \mathbb{R}^n\times\mathbb{R}^n\times \mathbb{R}.
  \end{equation*}
\end{description}
\end{theorem}
\bigskip

Let $p$ and $q$ be integers and
consider the simplest case where $n=1$,
$$V(x)=a_{2p}x^{2p}+a_{2p-1}x^{2p-1}+\cdots +a_1x+a_0,\ F(x)=c_{2(q+1)}x^{2(q+1)}+c_{2q+1}x^{2q+1}+\cdots +c_1x+c_0,$$
$E(x,y,t)\equiv e(t)\equiv e(t+\TT)$ and $\Sigma\equiv 0$. The necessary condition for system \eqref{DNewtPF} to be dissipative or the existence of a $\TT$-periodic solution is
\begin{equation}\label{leadingac}
a_{2p}>0, \qquad c_{2(q+1)}>0.
\end{equation}
Theorem \ref{LiWY} confirms such a condition is also sufficient when $p>q$.

To demonstrate the strength of Theorem \ref{UF}, we examine how far it can be applied to system (\ref{e:2.11}) with
$$F(x)=Q_{2q}(x)+\sum_{k=2}^{2q-1}Q_k(x), \qquad V(x)=P_{2p}(x)+\sum_{j=1}^{2p-1}P_j(x) , $$
where each $Q_k(x)$ and $P_j(x)$ is a homogeneous polynomial of degree $k$ and $j$, respectively.  In higher dimensions,  the scalar dissipative condition (\ref{leadingac})  corresponds to the requirement that
\vspace{1.0mm}\\
{\centerline{\it the leading terms $P_{2p}(x)$ and $Q_{2q}(x)$ are positive definite.}}
\vspace{0.001mm}\\
Under this minimal assumption alone,
we use Theorem \ref{UF} to prove the stochastic Levinson conjecture holds for polynomials $F$ and $V$, which greatly improves upon the counterpart result in \cite{LiWY}:
\medskip

\begin{theorem}\label{UF1}
Let $E(x,y,t), \Sigma(x,y,t)$ be continuous, $\TT$-periodic in time $t$,  locally Lipschitz continuous with respect to $(x,y)$, and $E(x,y,t)$ be bounded. Assume that
\begin{equation}\label{cC}
\| \Sigma\|_2^2(x,y,t) \le c\big(\|y\|^2+\|x\|^{2(p+q-1)}\big)+C  ,
 \end{equation}
 where $c$ is given in \eqref{cC1}. If both $Q_{2q}(x)$ and $P_{2p}(x)$ are positive definite and either $Q_{2q}(x)=\|x\|^{2q}$ or $P_{2p}(x)=\|x\|^{2p}$, then  system \eqref{e:2.11} admits a $\TT$-periodic probability solution.
\end{theorem}

  Note that the time-periodic Langevin equation (\ref{PLangrvin}) is included in general equation:
\begin{equation}\label{LEq}
\begin{cases}
   dx_t= y_tdt, \\
    dy_t=-\left[C(x_t,y_t,t)y_t +\nabla V(x_t)+E(x_t,y_t,t)\right] \, dt+\Sigma(x_t,y_t,t) \, d{B}_t.
\end{cases}
\end{equation}
Finally, we prove the existence of a periodic probability solution for (\ref{LEq}) with a uniformly positive friction and a potential of degree exceeding $2$.
\begin{theorem}\label{UF2}
The system \eqref{LEq} admits at least a $\TT$-periodic probability solution if it satisfies the Assumptions {\rm {\bf (A1)}-{\bf\bf(A4)}}. Here
\begin{description}
  \item{\rm \bf (A1)} All coefficients in system \eqref{LEq} are continuous in their respective domains and locally Lipschitz continuous in $(x,y)$. Let $V(x)\in C^{2}(\mathbb{R}^n)$ be lower bounded and $C(x,y,t),\ E(x,y,t)$,\ $\Sigma(x,y,t)$  also be periodic with respect to time $t$ with period $\TT$.

  \item{\rm \bf (A2)} There exist positive constants $\alpha$ and $\beta$ such that
 \begin{eqnarray*}
  C^s(x,y,t)\ge 2\alpha I_n;\ \ \|C(x,y,t)\|\, \vee \,|E(x,y,t)\| \le \beta, \qquad  \forall (x,y,t)\in \mathbb{R}^n\times\mathbb{R}^n\times \mathbb{R},
  \end{eqnarray*}
  where $C^s(x,y,t):=\frac{1}{2}\left(C(x,y,t)+C^{\top}(x,y,t)\right)$.

  \item{\rm \bf (A3)} There exist positive constants $b, M$ and $ \epsilon$ such that
  \begin{equation*}
  \langle x,\nabla V(x)\rangle \ge b\|x\|^{2+\epsilon}-M, \qquad\quad \forall x\in \mathbb{R}^n.
  \end{equation*}

  \item{\rm \bf (A4)} There exist positive constants $c$ and $M_1$ with $c<1$ such that
  \begin{equation*}
  \| \Sigma\|_2^2 (x,y,t)\le c \alpha\, (\|y\|^2 +b\|x\|^{2+\epsilon})+M_1, \qquad\qquad~ \forall (x,y,t)\in \mathbb{R}^n\times\mathbb{R}^n\times \mathbb{R}.
  \end{equation*}
\end{description}
\end{theorem}

We note that Theorems \ref{UF}-\ref{UF2} still hold under weaker regularity of the diffusion coefficient $A = (a^{ij}) \in C_{\TT}\left(\mathbb{R}, W_{{\rm loc}}^{1,p}(\mathbb{R}^n)\right)$ with all others continuous. We refer the reader to Theorem \ref{T:Qi} for details.

The proofs of Theorems \ref{UF}-\ref{UF2} rely on the periodic probability solution theory established by Ji {\color{red}et al.} \cite{Qi1,Qi2}  via a delicate Lyapunov function construction, building upon ideas from Jiang et al \footnote[1]{On the exponential ergodicity and large deviations of stochastic Hessian-driven damping Hamiltonian systems. Preprint.}. Our polynomial Lyapunov functions are smaller than the exponential functions introduced by Wu \cite{Wu}. Polynomial functions can offer superior probabilistic properties: their slower growth reduces the possibility of exponential explosion and broadens the admissible class of initial distributions.

A key contribution of this paper is the Lyapunov function defined in (\ref{e:3.2}), which combines mixed energy, potential energy, and a friction term. This formulation enables considerably more concise computations and provides the following advantages:
\smallskip

\begin{itemize}
\item {\bf  Stochastic Levinson's Conjecture Holds Under Broad Conditions}

\qquad All known results establishing the stochastic Levinson conjecture require the friction matrix to be uniformly positive definite, the potential force to diverge to  infinity and the noise to be bounded. Moreover, \cite{Liyong} and \cite{LiJ} additionally impose the assumption that the potential function $V(x)$
exhibits at most quadratic growth under condition (ii) of Theorem \ref{LiWY}. The result in \cite{Zhao} specifically addresses the case where the friction matrix is a positive constant multiple of the identity matrix. Works such as \cite{Qi1, Qi2} accommodate a broader class of friction matrices that are unbounded  in $x$ and bounded in $y$, but extracting the unbounded part $\Phi(x)$ of the friction matrix $C(x,y,t)$ required for their third assumption remains challenging.

\qquad
In contrast, our main Theorem \ref{UF} proves that distributed periodic solutions exist if $V(x,t)+aF(x) \to \infty$ uniformly in $t$ as $\|x\| \to \infty$ and the friction-potential gradient inner product grows at least like an even-power polynomial. The result holds even under bounded  friction matrix perturbations (requiring higher potential growth). This result substantially advances the stochastic Levinson conjecture, including cases not covered previously:
\begin{enumerate}
  \item stochastic relaxation oscillator: alternating friction definiteness;

  \item the Newton systems with friction matrix no lower bound;

  \item unrestricted growth rate in the friction-potential relationship;

  \item the Newton systems with bounded potential force;

  \item unbounded noises, potentially;

  \item bounded perturbations of Hessian friction: avoiding unbounded extraction.
\end{enumerate}
\end{itemize}

\begin{itemize}
\item {\bf Distributed Periodic Solutions Always Exist for Periodically Forced and Noisy Polynomial Systems}

\qquad Theorem \ref{UF1} establishes the existence of distributed periodic solutions for periodically forced and noise-driven polynomial systems with positive definite leading terms, which is a minimal assumption because it is the necessary condition for dissipation.

\qquad Theorem \ref{LiWY} requires strong dissipativity: the second inequality in (\ref{e:H02}) implies that the potential grows at least quadratically at infinity and has a higher degree than the friction (i.e., $p > q$). Consequently, it applies only to {\it strongly dissipative deterministic} time-periodic Newtonian systems. This excludes seminal cases like the time-periodically forced van der Pol equation \textemdash~a milestone in dynamical systems (see \cite{CL, Levinson1, Smale}). Furthermore, the first inequality in (\ref{e:H02}) mandates that friction admits a positive lower bound at a neighborhood of infinity.

\qquad By contrast, Theorem \ref{UF1} applies to {\it both deterministic and stochastic} time-periodic Newtonian systems, with no restrictions on the relative degrees of friction and potential. This broader applicability significantly extends Theorem \ref{LiWY}, even in the scalar case (see \cite[P. 2619-2620]{LiWY}). Specifically, it covers periodically forced and noise-driven van der Pol oscillator, van der Pol-Duffing oscillator, and their higher-dimensional counterparts. It even applies to dissipative systems with friction lacking a lower bound.

\qquad In particular, Theorem \ref{UF1} resolves Open Problem \eqref{V4} posed in \cite[Section 4]{Liyong}.
\end{itemize}
\begin{itemize}
\item {\bf Relaxing Infinity Requirement of Potential in Plasma Physics}

\qquad  Luigi Nocera  (MR4013829) (see his review on \cite{Qi1}) observed that existing results require the potential force to be infinite at infinity, which excludes certain plasma physics models. He suggested that distributed periodic solutions should exist for these forces and that the requirement could likely be relaxed for relevant plasma models.

\qquad Theorem \ref{UF} achieves precisely this relaxation: even when $\nabla_x V(x,t)$ is bounded, condition {\bf (H3)} holds provided $\nabla F(x)$ is so large that their gradient inner product grows like still an even-power polynomial, as demonstrated in \eqref{e:Exam1}-\eqref{e:Exam3} of Examples \ref{ex:4.2}--\ref{ex:4.4}. This further demonstrates the power of our Lyapunov function approach.
\end{itemize}

The organization of this paper is as follows.
In Section \ref{S:2}, necessary preliminaries
 are given and discussions on relations of various distributed periodic solution concepts are provided. In Section \ref{S:3}, we will give the proofs of Theorems \ref{UF}-\ref{UF2}.  In Section \ref{S:4}, we will provide a series of applications originating from plasma physics. In Section \ref{S:5},  we summarize our main results and state our further perspectives on stochastic time-recurrent Newtonian systems.

\section{Preliminaries}\label{S:2}

This section is devoted to provide preliminaries for our setup. Firstly, the definition of periodic probability solution or periodic solution in distribution is  introduced, which is referred to \cite{Liu1, LiJ, Liyong, Khasminskii, Qi1, Qi2, Zhao, Lu1, Lu2}. 
 Then, we present \cite[Theorem B]{Qi1},  which provides conditions for the existence
 of periodic probability solutions of SDEs in terms of Lyapunov functions.
\medskip

At first, we introduce the following notation.
Let
$(\Omega, \FF, \{\FF_t\}_{t \in \R}, \mathbb{P})$
be a filtered probability space satisfying the usual conditions. 
For periodic stochastic systems,
the period of function $f$ with respect to time $t$ being $\TT$ means that
$
f(x,t+\TT) =f(x,t)
$
holds for all $t\in \R$ and $x\in \R^n$.
Denote by $ C_\TT ( \R^n\times  \R)$ the space of $\TT$-periodic and continuous functions on $ \R^n \times \R$ and by $C_\TT^{2,1} ( \R^n\times \R)$ the space $ C^{2,1} ( \R^n\times  \R) \cap C_\TT ( \R^n\times  \R)$. $C_0^{2,1} ( \R^n\times \R)$ denotes the space of all those functions in $C^{2,1} ( \R^n\times \R)$ with compact supports.
\smallskip

Ji et al \cite{Qi1, Qi2} considered the following $\TT$-periodic stochastic differential equation
\begin{equation}\label{e:2.4}
\begin{cases}
  dZ_t= b(Z_t,t)\, dt+ \sigma (Z_t,t)\, dB_t, ~~\quad t>s ,\\
  Z_s= z
\end{cases}
\end{equation}
with Sobolev coefficients. They studied periodic solutions to (\ref{e:2.4}) with irregular (in particular, non-Lipschitz) coefficients in the sense of distribution by its corresponding
 Fokker-Planck equation:
 \begin{equation}\label{FP}
\mathcal{L}^*p := -\partial_tp + \partial_{ij}^2\left(a^{ij}p\right) - \partial_i(b^ip)=0, (x, t) \in \mathbb{R}^n\times \mathbb{R}.
\end{equation}
Here the diffusion matrix $(a^{ij})=\frac12 \sigma \sigma^T$ and the drift field $b$ are $\TT$-periodic in its second variable,  $\partial_i=\partial_{z_i}, \partial_{ij}=\partial_{z_iz_j}, i, j\in \{1, 2, \cdots, n\}$,
and the usual summation convention is used. The following definitions are adopted from \cite[p.4-5]{Qi1}.
\begin{definition}{\rm(Periodic probability solution)}\label{PPS}
A Borel measure $\mu$ on $\mathbb{R}^n\times \mathbb{R}$ is called a
{\it periodic probability solution} to {\rm(\ref{FP})} if there is a family of Borel probability measures $\{\mu_t\}_{t\in \mathbb{R}}$ on $\mathbb{R}^n$ satisfying
$$\mu_t=\mu_{t+\TT},\ \forall t\in \mathbb{R},$$
$$a^{ij}, b^i\in L_{{\rm loc}}^1\left(\mathbb{R}^n\times \mathbb{R}, {\rm d}\mu_t{\rm d}t\right),\ \forall i, j\in \{1, 2, \cdots, n\} $$
and $\mathcal{L}^*\mu=0$ on $\mathbb{R}^n\times \mathbb{R}$ in the sense that
\begin{equation}\label{PMS}
\int_{\mathbb{R}}\int_{\mathbb{R}^n}\mathcal{L}\phi{\rm d}\mu_t{\rm d}t=0,\ \forall \phi\in C^{2,1}_0 ( \R^n\times  \R)
\end{equation}
such that ${\rm d}\mu={\rm d}\mu_t{\rm d}t$.
\end{definition}

For  a non-negative function $\Psi\in C_\TT^{2,1} ( \R^n\times \R)$ and each $\rho>0$, we define the $\rho-$sublevel set
$$\Omega_{\rho}:=\{(z,t)\in \mathbb{R}^n\times \mathbb{R}: \Psi(z,t)<\rho\}.$$
\begin{definition}{\rm(Unbounded Lyapunov function)}\label{ULF}
A non-negative function $\Psi\in C_\TT^{2,1} ( \R^n\times \R)$ is called an unbounded Lyapunov function with respect to $\mathcal{L}$ if there is a sequence $ \{\mathcal{U}_n\}_{n\in\mathbb{N}}$ of open sets in $\mathbb{R}^n$ satisfying $\mathcal{U}_m\subset \mathcal{U}_{m+1}\subset\subset \mathbb{R}^n$ for all $m\in \mathbb{N}$ and  $\mathbb{R}^n=\bigcup_{m=1}^{+\infty}\mathcal{U}_m$
such that
\begin{equation}\label{Pinfty}
{\rm inf}_{(\mathbb{R}^n\backslash \mathcal{U}_m)\times \mathbb{R}}\Psi(z,t)\rightarrow +\infty,\ {\rm as}\ m\rightarrow +\infty,
\end{equation}
and there exist a $\rho_k>0$, called an essential lower bound of $\mathbb{R}^n$, and a constant $\gamma>0$,
called a Lyapunov constant of $\Psi$, such that
\begin{equation}\label{Ninfty}
\mathcal{L}\Psi(z,t)\le -\gamma,\ \forall (z,t)\in \left(\mathbb{R}^n\times \mathbb{R}\right)\backslash {\Omega_{\rho_k}}.
\end{equation}
\end{definition}

The following is  the existence theorem
of periodic probability solutions to the Fokker-Planck equation (\ref{FP}) with Sobolev coefficients, due to \cite{Qi1}.
\begin{theorem}\label{T:Qi}
{\rm(\cite[Theorem B]{Qi1})}
Let {\bf (H)} hold. Suppose $A = (a^{ij}) \in C_{\TT}\left(\mathbb{R}, W_{{\rm loc}}^{1,p}(\mathbb{R}^n)\right)$ and $b\in C_{\TT}\left(\mathbb{R}^n\times \mathbb{R}\right)$. If there is an unbounded Lyapunov function with respect to $\mathcal{L}$, then {\rm(\ref{FP})} admits a periodic probability solution. Here the hypothesis {\bf (H)} is given as follows.

{\bf (H)} Let $p > n + 2$ and $ a^{ij} \in  L^{\infty}_{{\rm loc}}\left(\mathbb{R}, W_{{\rm loc}}^{1,p}(\mathbb{R}^n)\right)
$  for each $ i, j \in \{1, ..., d\}$. The drift
vector field $b= (b^i)$ satisfies $b^i\in L^p_{\rm loc}\left(\mathbb{R}^n\times \mathbb{R}\right)$ for each $ i\in \{1, ..., d\}$.
\end{theorem}

\begin{cor}\label{T:Khasminskii}
Suppose that the coefficients of \eqref{e:2.4} are $\TT$-periodic in time $t$ and satisfy all regular assumptions in Theorem {\rm\ref{T:Qi}}. If there exists a non-negative
function $ \Psi(x,t) \in C_\TT^{2,1} ( \R^n, \R)$, which satisfies
\begin{equation}\label{condition2}
\lim_{\|z\|\to \infty} ~  \Psi(z,t) = +\infty,\ {\rm uniformly\ in}\ t\in \mathbb{R}
\end{equation}
and
\begin{equation}\label{condition1}
\lim_{\|z\|\to \infty} ~  \mathcal{L} \Psi (z,t) = -\infty,\ {\rm uniformly\ in}\ t\in \mathbb{R},
\end{equation}
then {\rm(\ref{FP})} admits a periodic probability solution.
\end{cor}
Note that (\ref{condition2}) and (\ref{condition1}) obviously imply (\ref{Pinfty}) and (\ref{Ninfty}), respectively. Under the assumption of locally Lipschitz continuity of coefficients, Corollary \ref {T:Khasminskii} is \cite[Theorem 3.8]{Khasminskii}.

 From now on we assume that the drift $b: \R^n \times \R \to \R^n$ and noise intensity $ \sigma: \R^n \times \R \to \R^{n\times n}$ are continuous on $\R^n \times \R$, locally Lipschitz continuous in their first variable $z$ and  $\TT$-periodic in their second variable $t$. Under this assumption, the transition probabilities associated to solutions of \eqref{e:2.4} can be defined.  Let $ Z^{s,z}(t)$ be the solution of (\ref{e:2.4}) and $\P(s, z, t, A):= \mathbb{P}\big(Z^{s,z}(t)\in A\big) $ its transition probability function. Then $\P(s, z, t, A)$ is $\TT$-periodic in the sense that $\P(s+\TT, z, t+\TT, \cdot)\equiv\P(s, z, t, \cdot)$ for all $z\in \mathbb{R}^n$ and $s\le t$, see \cite[Theorem 3.4]{Khasminskii}.
We call a Markov process $ Z(t)$ to be {\it $\TT$-periodic} if its law $ \P_0 ( t, A):= \P\{ Z(t) \in A\}$ satisfies the equation
\begin{equation}\label{periodicS}
\P_0 ( s, A)= \int_{\R^n} \P_0 ( s, dz) \, \P(s, z, s+\TT, A)
\equiv \P_0 ( s+\TT, A), \quad  \forall s\in\mathbb{R},\ A\in \mathcal{B}(\mathbb{R}^n),
\end{equation}
see \cite[P. 64, \,(3.17)]{Khasminskii}. (\ref{periodicS}) is equivalent to that the laws of $Z(s)$ and $Z(s+\TT)$ are identical for any $s\in\mathbb{R}$. According to \cite{LiJ, Liyong}, such a $\TT$-periodic solution $Z(t)$ is called to be {\it periodic solution in distribution}.
\smallskip

Ji et al \cite[p.3-4]{Qi1} commented that the distribution of solutions to \eqref{e:2.4} is governed by (\ref{FP}) at least when the coefficients are locally Lipschitz continuous. Therefore, the distribution of periodic solutions to \eqref{e:2.4} correspond to periodic probability solutions to (\ref{FP}). We will give a proof for such a statement and further discuss the relations between the concepts of periodic probability solutions and periodic solutions in distribution.

Suppose that $(Z(t))_{t\in\mathbb{R}}$ is a solution of \eqref{e:2.4} whose law or distribution is $\TT$-periodic, i.e., let $\mu_t=\mathbb{P}\circ(Z(t))^{-1}$, then $\mu_t=\mu_{t+\TT},\ \forall t\in \mathbb{R}$. We claim that (\ref{PMS}) holds.

For a given $\phi\in C^{2,1}_0 ( \R^n\times  \R)$, we can find a positive integer $k$ such that ${\rm supp}(\phi)\subset \overline{B_k(O)}\times [-k\TT,k\TT]$.
An application to $\phi$ of It\^o's formula yields
\begin{equation*}
\begin{split}
\mathbb{E}_{\mu_s} \phi (Z(t),t)-\int_{\mathbb{R}^n} \phi(\cdot, s) {\rm d} \mu_s&=\int_s^t \mathbb{E}_{\mu_s} \mathcal{L}\phi (Z(r),r) {\rm d} r\\
&=\int_s^t \int_{\mathbb{R}^n} \mathcal{L}\phi(z,r) {\rm d} \mu_r{\rm d}r,\quad s\leq t.
\end{split}
\end{equation*}
Let $s<-k\TT$ and $t>k\TT$. It follows from ${\rm supp}(\phi)\subset \overline{B_k(O)}\times [-k\TT,k\TT]$ and the above equality that
$$
0=\int_s^t \int_{\mathbb{R}^n} \mathcal{L}\phi(z,r) {\rm d} \mu_r{\rm d}r=\int_{\mathbb{R}} \int_{\mathbb{R}^n} \mathcal{L}\phi(z,r) {\rm d} \mu_r{\rm d}r.
$$
This proves the claim, i.e., $(\mu_t)_{t\in\mathbb{R}}$ is a periodic probability solution of (\ref{FP}) according to \cite{Qi1, Qi2}. This indicates that their periodic probability solution theory generalizes the classical periodic solutions of SDEs and accommodates coefficients with merely Sobolev regularity.

 Conversely, assume that $(\mu_t)_{t\in\mathbb{R}}$ is a periodic probability solution of (\ref{FP}) according to \cite{Qi1}, let $ Z^{s}(t),\ s\le t$ be a solution of \eqref{e:2.4} with initial distribution $\mu_s$ and $\tilde{\mu}^s_t$ the distribution of $Z^s(t)$. In particular, $\tilde{\mu}^s_s=\mu_s$. Then we have
\begin{equation}\label{periodicS1}
\tilde{\mu}_{s+\TT}(A)= \int_{\R^n} \mu_s(dz) \, \P(s, z, s+\TT, A),\quad  \forall \ A\in \mathcal{B}(\mathbb{R}^n).
\end{equation}
An application of It\^o's formula yields
$$
\mathbb{E}_{\mu_s} \phi (Z^s(t),t)-\int_{\mathbb{R}^n} \phi(\cdot,s) {\rm d} \mu_s=\int_s^t \mathbb{E}_{\mu_s} \mathcal{L}\phi (Z^s(r),r ) {\rm d} r=\int_s^t \int_{\mathbb{R}^n} \mathcal{L}\phi {\rm d} \tilde{\mu}^s_r{\rm d}r,\quad \forall\ \phi\in C^{2,1}(\mathbb{R}^d\times \mathbb{R}).
$$
In other words,
\begin{equation}\label{NEFP}
 \int_{\mathbb{R}^n}\phi (\cdot,t){\rm d} \tilde{\mu}^s_t=\int_{\mathbb{R}^n} \phi(\cdot,s) {\rm d} \tilde{\mu}^s_s + \int_s^t \int_{\mathbb{R}^n} \mathcal{L}\phi {\rm d} \tilde{\mu}^s_r{\rm d}r,\quad \forall\ s\le t,\ \phi\in C^{2,1}(\mathbb{R}^d\times \mathbb{R}).
\end{equation}
That is, $\tilde{\mu}^s_t$ satisfies \eqref{FP} with the initial distribution $\mu_s$ in the sense of distribution. Note that \cite[Lemma 2.1 (1)]{Qi2} ensures that $\{\tilde{\mu}^s_t\}_{t\geq s}$  is also a global probability solution of \eqref{FP} with the  initial distribution $\mu_s$. Suppose that the Cauchy problem associated with \eqref{FP} admits a unique global probability solution. Then $\mu_t=\tilde{\mu}^s_t=\mathbb{P}[Z^s(t)\in \bullet]$ for all $t\geq s$. This shows that
\begin{equation*}\label{periodicS1}
\mu_{s+\TT}(A)= \int_{\R^n} \mu_s(dz) \, \P(s, z, s+\TT, A)=\mu_{s}(A),\quad  \forall\ s\in\mathbb{R}, \ A\in \mathcal{B}(\mathbb{R}^n).
\end{equation*}
This illustrates the following: under the assumptions that the coefficients are locally Lipschitz continuous and that the probabilistic solution for Cauchy problem of system (\ref{FP}) possesses global uniqueness, the concept of a periodic probability solution in \cite{Qi1,Qi2} coincides with the concept of a periodic solution in distribution in \cite{LiJ, Liyong}.
Under the non-degenerate assumption, \cite[Corollary A]{Qi2} guarantees that such a uniqueness of global probability solution. However, stochastic Newtionian systems (\ref{LEq}) are always degenerate.  Fortunately, in many cases, it can be verified that (\ref{NEFP}) satisfies the H\"{o}rmander condition so that the Fokker-Planck operator becomes hypoelliptic. Given the similarities between hypoelliptic and elliptic operators, it seems promising to prove the uniqueness of global probability solutions, under additional smooth assumption of coefficients and Lyapunov dissipative condition. We leave it for future study.

Based on the above discussion, we adopt the more broadly applicable {\it periodic probability solution} as the definition for distributed periodic solution throught this paper.

The following example  is taken from \cite[subsection 3.1]{Liyong}. They made many numerical experiments about it. Their results shows the solutions starting from (0, 0) do not admit periodicity pathwisely. However, the sample paths do admit periodic phenomenon significantly in distribution. These phenomena confirm them that it is
natural and essential to study stochastic dynamical problems in the sense of distribution. In the following, we explicitly solve its density of the $2\pi-$periodic probability solution or $2\pi-$periodic solution  in distribution $(\mu_t)_{t\in\mathbb{R}}$ for \eqref{Duffing}, which helps readers to understand these concepts.

\begin{example}\label{lizi} \rm
Consider the scalar periodic stochastic Duffing equation:
$$
\ddot{x} + \dot{x} + 2x=\dot{B}_t+ \cos t.
$$
Equivalently, we have
\begin{eqnarray}\label{Duffing}
\begin{cases}
   dx_t= y_t \,dt,  \\
    dy_t=(-2x_t-y_t+\cos t)\,dt+ dB_t.
\end{cases}
\end{eqnarray}
Then define $Z:= (x, y)^\top$, it can be rewritten as
$$
dZ_t = ( AZ_t +b(t)) \,dt+ \Sigma \,dB_t ,
$$
where
\begin{eqnarray*}
A= \left(\begin{matrix}
 0 &1 \\
 -2 &-1
\end{matrix} \right),  \quad
 b(t)=\left(\begin{matrix}
 0\\
\cos t
\end{matrix}
\right) ,\quad
\Sigma= \left(\begin{matrix}
 0\\
1
\end{matrix}
\right)  \quad
~\hbox{and}\quad~ \Sigma\Sigma^\top
=\left(\begin{matrix}
 0 &0 \\
 0 &1
\end{matrix} \right).
\end{eqnarray*}
\end{example}
Its  Fokker-Planck equation is
\begin{eqnarray}\label{FP1}
\frac{\partial p(x,y,t)}{\partial t}=\frac12 \frac{\partial^2 p(x, y,t)}{\partial {y}^2} - y \frac{\partial p(x,y,t)}{\partial x} + (2x+ y-\cos t ) \frac{\partial p (x, y,t)}{\partial y } + p(x, y,t).
\end{eqnarray}
Explicitly, we can obtain
\begin{eqnarray*}
  p(x,y,t)= \frac{\sqrt 2}{\pi}\, \exp{\left(-2 \left(x- \frac{\cos t+\sin t }{2} \right)^2- \left( y -
  \frac{\cos t-\sin t }{2}  \right)^2 \right)},
\end{eqnarray*}
 which is the density function  of  $2\pi$-periodic probability solution $(\mu_t)_{t\in\mathbb{R}}$ for system \eqref{Duffing}.
\qed
\medskip

\section{The Proof of Main Theorems}\label{S:3}

{\bf Proof of Theorem \ref{UF}.} \quad
Define the Lyapunov function
\begin{eqnarray}\label{e:3.2}
\Psi(x,y,t):= \frac12 \| y+ \nabla F(x) -ax\|^2  + \Big[ V(x,t) +aF(x) -\frac{a^2}{2} \| x  \|^2  \Big] +D.
\end{eqnarray}
Here, $a>0$ is given in Hypothesis {\bf(H2)} and related to convergent rate, and $D$ is a sufficiently large constant such that $ \Psi\geq 1$ for any $ (x,y,t) \in \R^n \times \R^n \times \R$. The first term is called to be the {\it mixed energy}. The selection of cross term between variables $x$ and $y$ is absolutely crucial when constructing Lyapunov functions for damped Newtonian systems. We choose the cross term to be $\langle y, \nabla F(x) - a x \rangle$ within the mixed energy, leading to significantly more concise computations and greatly advancing  the applicable scope as described in Introduction.

Obviously, function $\Psi$ is $\TT$-periodic and satisfies by Hypothesis {\bf(H2)}
$$
\lim_{\|x\|+\|y\|\rightarrow +\infty} \Psi(x,y,t)= +\infty
$$
uniformly in time $t\in \mathbb{R}$.
This implies function $\Psi$ satisfies condition \eqref{condition2} of Corollary \ref{T:Khasminskii}.

The generator of system \eqref{e:2.11} is given by: for any $f\in C^{1,2,1}: \R^n \times\R^n \times \R \to \R$,
\begin{eqnarray}
&&\mathcal{L} f(x,y,t) \nonumber\\
\label{e:3.3}
&=& \partial_t f(x,y,t)+\frac12 \sum_{i,j=1}^n (\Sigma  \Sigma^\top)^{i,j} (x,y,t) \, \partial^2_{y_i\,y_j } 
 f(x,y,t)  \nonumber\\
&& + \<y, \nabla_x f(x,y,t) \> - \<D^2 F(x)\, y+ \nabla_x V(x,t)+E(x,y,t) , \nabla_y f(x,y,t) \>  .
\end{eqnarray}
It remains to check condition \eqref{condition1} of Corollary \ref{T:Khasminskii}. By computation, we have
\begin{eqnarray}\label{y}
\nabla_y\Psi=y +\nabla F(x) -ax,
\end{eqnarray}
\begin{eqnarray}\label{x}
\nabla_x\Psi &=& [\nabla_xV(x,t)+a\nabla F(x)-a^2x]\nonumber\\
&&+ [D^2F(x)y-ay+D^2F(x)\nabla F(x)-a\nabla F(x)-aD^2F(x)x+a^2x]\nonumber\\
&=& \nabla_xV(x,t)+ D^2F(x)\,y-ay+D^2F(x)\nabla F(x)-aD^2F(x) \,x.
\end{eqnarray}
Then using \eqref{y} and \eqref{x}, we get
\begin{eqnarray}
&&\<y, \nabla_x \Psi(x,y,t) \> - \<D^2 F(x)\, y+ \nabla_x V(x,t) , \nabla_y \Psi(x,y,t) \>
 \nonumber\\
 \label{keyterm}
&=&-\big[a\|y\|^2+\<\nabla_x V(x,t), \nabla F(x)-ax\>\big].
\end{eqnarray}
Substituting \eqref{keyterm} into \eqref{e:3.3} with $f=\Psi$, and then applying Hypotheses {\bf(H3)}-{\bf(H5)}, we obtain that
\begin{eqnarray}
&&\mathcal{L} \Psi(x,y,t)\nonumber \\
\label{L3}
&=& \partial_t V(x,t) +\frac12 \| \Sigma\|_2^2 (x,y,t)  -\<E(x,y,t), y+\nabla F(x)-ax  \>
 \nonumber\\
&&-a\|y \|^2 - \<\nabla_x V(x,t), \nabla F(x)-ax \> \nonumber  \\
&\leq& (c_1 +c_2 -1)\,( a\|y\|^2 + b \|x\|^{2m}) + \Big(M_1+ \frac12 M_2+M \Big) \nonumber \\
&& +e\|y\| +e\|\nabla F(x)-ax \| \label{L3}  \\
&\leq& \frac{1}{2}(c_1 +c_2 -1)\,( a\|y\|^2 + b \|x\|^{2m}) + \Big(M_1+ \frac12 M_2+M +M^*\Big)
\label{L4}
\end{eqnarray}
for some positive constant $M^*$, where the first inequality \eqref{L3} has used Hypotheses ({\bf H5}), ({\bf H3}) and $\|E(x,y,t)\|\le e$, while the second inequality \eqref{L4} has used Hypothesis ({\bf H4}). This demonstrates \eqref{condition1} from the fact that $0<c_1 +c_2 <1$ and then the conclusion follows immediately from Corollary \ref{T:Khasminskii}.
\qed
\vspace{2.0ex}\\
{\bf Proof of Theorem \ref{UF1}.} \quad The proof of this theorem is mainly obtained by applying Theorem \ref{UF}, so we only need to verify Hypotheses {\bf(H1)}-{\bf(H5)} one by one.  The proof is divided into several steps. It is easy to see that Hypothesis {\bf(H1)} is satisfied by our assumptions here.
\bigskip

{\bf Step 1.} We first prove that Hypothesis {\bf(H2)} holds.
In order to check it, we have to find a positive constant $a$ such that
\begin{equation}\label{Pinfinity1}
  \lim_{\|x\|\rightarrow +\infty} \Big[V(x)+aF(x)-\frac{a^2}{2}\|x\|^2 \Big]=+\infty.
  \end{equation}

Indeed, if either $p>1$ or $q>1$, then by the positive definiteness of $Q_{2q}(x)$ and $P_{2p}(x)$, \eqref{Pinfinity1} holds for any given positive $a$. It remains to consider the case $p=q=1$. In this case, $Q_{2}(x)$ and $P_{2}(x)$ are positive definite quadratic forms. Thus, there exists a positive constant $\lambda$ such that
\begin{eqnarray*}
Q_{2}(x)\, \wedge \, P_{2}(x)\ge \lambda \|x\|^2.
\end{eqnarray*}
If necessary, we can shrink $\lambda$ so that
$$\Big[V(x)+aF(x)-\frac{a^2}{2}\|x\|^2 \Big]\ge \Big(\lambda+\lambda a -\frac{a^2}{2} \Big)\, \|x\|^2 $$
holds. Choose $a\in (0, 2\lambda)$ so small that $(\lambda+\lambda a -\frac{a^2}{2})>0$. Such an $a$ guarantees that \eqref{Pinfinity1} holds. Therefore, there always exists a positive constant $a$ such that \eqref{Pinfinity1} holds, which proves that Hypothesis {\bf(H2)} is always true.
\medskip

{\bf Step 2.} We then claim that there exist positive constants $b$ and $M$ such that
  \begin{equation}\label{H3}
  \langle \nabla V(x), \nabla F(x)-ax\rangle \ge b\|x\|^{2(p+q-1)} - M, \qquad \forall x\in \mathbb{R}^n.
  \end{equation}
This implies that Hypothesis {\bf(H3)} holds.

  Without loss of generality, we assume that  $P_{2p}(x)=\|x\|^{2p}$.   First, let $p=q=1$. Then by the Euler formula, we can obtain
  \begin{eqnarray*}
&&\langle \nabla V(x), \nabla F(x)-ax\rangle \\
&=& \langle 2x+v_0, \nabla Q_2(x)-ax\rangle\\
&=& 4Q_2(x)-2a\|x\|^{2}+\langle v_0, \nabla Q_2(x) \rangle-a\langle v_0,x \rangle  \\
&\geq& 2(2\lambda-a)\|x\|^{2}+\langle v_0, \nabla Q_2(x) \rangle-a\langle v_0,x \rangle\\
&\geq& (2\lambda-a)\|x\|^{2} - M
\end{eqnarray*}
for some positive constant $ M$, which implies (\ref{H3}) holds in this case.

In the remaining case $p+q>2$, it follows from the positive definiteness of $Q_{2q}(x)$ that there is a positive constant $\nu$ such that $Q_{2q}(x)\ge \nu \|x\|^{2q}$ for all $x\in \R^n$. Employing the Euler formula, we deduce that
$$\langle x,\nabla Q_{2q}(x) \rangle =2q \, Q_{2q}(x)\ge 2q\nu\|x\|^{2q},
$$
and further
\begin{eqnarray*}
&&\langle \nabla V(x), \nabla F(x)-ax\rangle \\
&=& \big\langle 2p\|x\|^{2(p-1)}x+\sum_{j=1}^{2p-1}\nabla P_j(x),\nabla Q_{2q}(x)+\sum_{k=2}^{2q-1}\nabla Q_k(x)-ax \big \rangle\\
&=& 2p\, \|x\|^{2(p-1)}\langle x, \nabla Q_{2q}(x)\rangle+ \hbox{ lower-order\ terms}\\
&=& 4pq \, \|x\|^{2(p-1)} Q_{2q}(x)+\hbox{ lower-order\ terms}\\
&\ge& 4pq \, \nu\|x\|^{2(p+q-1)}+o(\|x\|^{2(p+q-1)})\ \quad ({\rm as}\ \|x\|\rightarrow +\infty)\\
&\ge& 2pq \, \nu\|x\|^{2(p+q-1)} - M.
\end{eqnarray*}
This proves \eqref{H3} for the case $p+q>2$.
\medskip

{\bf Step 3}. Finally, we verify Hypotheses {\bf\bf(H4)} and {\bf\bf(H5)}.
Note that a positive constant $c$ in \eqref{cC} should satisfy
\begin{equation}\label{cC1}
c<
\begin{cases}
   {\rm min}\{2(2\lambda-a), 2a\},\ {\rm if}\ p=q=1, \\
   {\rm min}\{4pq\nu, 2a\},\ {\rm if}\ p+q>2.
\end{cases}
\end{equation}
Thus far, we have verified Hypotheses {\bf (H1)}-{\bf\bf(H4)}. It is straightforward to verify that Hypothesis {\bf\bf(H5)} holds for
$m=2(p+q-1)$.
\smallskip

Based on the above displays, the conclusion of this theorem follows immediately from Theorem \ref{UF}.
\qed

\begin{remark}
In the deterministic situation,
Li et al {\rm\cite {LiWY}} considered
\begin{eqnarray}\label{DFV}
\ddot{x} + D^2F(x)\dot{x} + \nabla V(x)=e(t)\equiv e(t+\TT) ,
\end{eqnarray}
and gave sufficient conditions for system \eqref{DFV} to admit a $\TT$-periodic solution. For scalar case, that is, $n=1$, and the friction and the potential are polynomials:
$$V(x)=a_0+a_1x+a_2x^2+\cdots a_{2p}x^{2p}(a_{2p}>0),\ F(x)=c_0+c_1x+c_2x^2+\cdots c_{2(q+1)}x^{2(q+1)}(c_{2(q+1)}>0),$$
they proved that system \eqref{DFV} has a $\TT$-periodic solution if $p>q\ge 0$. For higher dimensional system, their result applies to $F(x)=\|x\|^{2q},\ V(x)=\|x\|^{2p}$ to conclude that the existence of  a $\TT$-periodic solution if $p>q>0$. However, it is worth noting that their result cannot be applied to the well-known periodically forced  van der Pol equation and its higher-dimensional version.

Note that system \eqref{DFV} with $n=1$ is dissipative if and only if  $a_{2p}>0$ and $ c_{2(q+1)}>0$.
For higher-dimensional systems, this condition corresponds to the leading terms of $V(x)$ and $F(x)$ being positive definite. Under this minimal assumption, we proved in Theorem {\rm\ref{UF1}} that both systems \eqref{DFV} and \eqref{e:2.11} admit a $\TT$-periodic solution. This result significantly improves upon those in {\rm\cite{LiWY}}.
\end{remark}
\vspace{1.5ex}
{\bf Proof of Theorem \ref{UF2}.}
 Setting $F(x):=\alpha \|x\|^2,\ a:=\alpha$ in \eqref{e:3.2} leads to the following Lyapunov function
$$\Psi(x,y)=\frac{1}{2}\|y+\alpha x\|^2+V(x)+\frac{\alpha^2}{2}\| x\|^2.$$
By Hypothesis ({\bf A1}), it is easy to see that
$$\lim_{\|x\|+\|y\|\rightarrow +\infty} \Psi(x,y)=+\infty.$$
This implies the function $\Psi$ satisfies condition \eqref{condition2} of Corollary \ref{T:Khasminskii}.
 In addition,  using Hypotheses {\bf (A2)}-{\bf (A4)}, we have
\begin{align*}\label{LV}
\mathcal{L}\Psi = &-\alpha [ \|y\|^2  + \langle x,\nabla V\rangle (x)] -\langle y,(C^s-2\alpha I)y\rangle  -\alpha \langle x,(C-2\alpha I)y\rangle -\langle E,y+\alpha x\rangle+\frac{1}{2}\|\Sigma\|_2^2\\
\le & -\alpha \big[\|y\|^2+\langle x,\nabla V(x)\rangle \big]+ \alpha (\beta+2\alpha)\|x\| \, \|y\|+ \beta(\|y\|+\alpha \|x\|) + \frac{1}{2}\|\Sigma\|_2^2\\
\le & -\alpha \big[\|y\|^2+\langle x,\nabla V(x)\rangle \big]+ \frac{\alpha}{4}\|y\|^2+\alpha (\beta+2\alpha)^2\|x\|^2+ \beta(\|y\|+\alpha \|x\|) + \frac{1}{2}\|\Sigma\|_2^2\\
 \le & -\frac{3\alpha}{4} \|y\|^2 +\beta \|y\|-\alpha b\|x\|^{2+\epsilon} + \alpha (\beta+2\alpha)^2\|x\|^2+ \alpha \beta \|x\| +\alpha M+\frac{1}{2}\|\Sigma\|_2^2\\
 \le & -\frac{\alpha}{2}\left( \|y\|^2 + b\|x\|^{2+\epsilon}\right) + M^*+\frac{1}{2}\|\Sigma\|_2^2\\
\le & -\frac{\alpha}{2}(1-c)\left( \|y\|^2 + b\|x\|^{2+\epsilon}\right) + M^*+\frac{M_1}{2},
\end{align*}
which implies the function $\Psi$ satisfies condition \eqref{condition1}. The conclusion follows immediately from Corollary \ref{T:Khasminskii}.
\qed

\begin{remark}
We note that Ji et al {\rm\cite{Qi1, Qi2}}  have established Theorem {\rm\ref{UF2}} under bounded noise conditions when the friction matrix is bounded, using the exponential-type Lyapunov function introduced in {\rm\cite{Wu}}. However, our approach employs a polynomial Lyapunov function. A polynomial function offers superior probabilistic properties compared to an exponential one. Specifically, its slower growth rate mitigates the risk of exponential explosion and expands the admissible space of initial distributions.
\end{remark}

\section{Plasma Physics Models with Bounded Potential Force}\label{S:4}
This section is motivated by Luigi Nocera's remarks on \cite{Qi1}, see MR4013829. He observed that existing results require potentials to be ``infinite at infinity,'' which excludes applications to certain plasma physics models. Nocera suggested that distributed periodic solutions should exist for these models and that the requirement could likely be relaxed. Our Theorem \ref{UF} achieves precisely this relaxation.

We present a series of examples with potential functions drawn from \cite{SSK, KMT, RSC} (sometimes slightly modified). In these cases, the potential force is bounded, while friction and noise terms remain unbounded, which further highlights the strength of our theorem. Throughout, we assume the diffusion coefficient $\Sigma$ satisfies the standard assumptions in Hypothesis {\bf(H1)}.
\smallskip

\begin{example}\label{ex:4.2}
\rm Suppose that
\begin{eqnarray*}
V(x,t)= \ln (2+\sin t + \|x\|^2 ) , \qquad F(x)= \| x\|^4 +\|x\|^2.
\end{eqnarray*}
We claim that if $  \| \Sigma\|_2^2 \le c_3 (\|x\|^2 +\|y\|^2 )  +C$ with $c_3 <16$, then it admits a $2\pi$-periodic probability solution.
\end{example}

In fact,
taking $a:=8$ and $ D:=19$,  we have
\begin{eqnarray*}
\partial_t V(x,t) &=& \frac{\cos t}{ 2+\sin t + \|x\|^2}, \\
\nabla_x V(x,t) &=&  \frac{2x}{ 2+\sin t + \|x\|^2},  \\
\nabla F(x)- 8x &=& 2x\, (2\|x\|^2 -3).
\end{eqnarray*}
Obviously, $\|\nabla_x V(x,t) \| \leq 1$ and $D^2 F(x)=(4\|x\|^2 +2) I_n + 8xx^{\top}$ is positive definite. Moreover,
\begin{equation}\label{e:Exam1}
-\< \nabla_x V(x,t), \nabla F(x)- 8x \> = - \frac{4\|x\|^2 \, (2\|x\|^2 -3 )}{2+\sin t + \|x\|^2} ~\leq -8\|x\|^2+36.
\end{equation}
It is easy to see that Hypotheses {\bf(H1)}-{\bf(H5)} hold with $ b=8, \, m=1, $ and $ M=36$.
\qed

\begin{example}\label{ex:4.3}
\rm Assume that
$$ V(x,t)= \sqrt{2+\sin t + \|x\|^2 } , \qquad F(x)= \| x\|^4 +\|x\|^2. $$
 We claim that if $  \| \Sigma\|_2^2 \le c_4 (\|x\|^3 +\|y\|^2 )  +C$ with $c_4 <4$, then it admits a $2\pi$-periodic probability solution.
\end{example}

Indeed, taking $a:=2$ and $ D:=1$, we obtain
\begin{eqnarray*}
\partial_t V(x,t) &=& \frac{\cos t}{ 2\sqrt{2+\sin t + \|x\|^2 } }, \\
\nabla_x V(x,t) &=&  \frac{x}{ \sqrt{2+\sin t + \|x\|^2 } },  \\
\nabla F(x)- 2x &=& 4 \|x\|^2\,x.
\end{eqnarray*}
Meanwhile,
$\|\nabla_x V(x,t) \| \leq 1$ and $D^2 F(x)=(4\|x\|^2 +2) I_n + 8xx^{\top}$ is positive definite.
In addition, it yields
\begin{equation}\label{e:Exam2}
-\< \nabla_x V(x,t), \nabla F(x)- 2x \> = - \frac{ 4 \|x\|^4}{ \sqrt{2+\sin t + \|x\|^2 } } \leq -2\|x\|^3 +2 .
\end{equation}
Hence, Hypotheses {\bf(H1)}-{\bf(H5)} hold with $ b=2, \, m=3/2, $ and $ M=2$.
Then we conclude from Theorem \ref{UF} that there exists a $2\pi$-periodic probability solution.
\qed

\begin{example}\label{ex:4.4}
\rm Let
$$
V(x,t)= (2+\sin t)(1- \exp{(-\|x\|^2) } ), \qquad F(x)= \frac14 \int_0^{ \|x\|^2} e^s\, (s+1)\, ds.
$$
We claim that if $  \| \Sigma\|_2^2 \le c_5 (\|x\|^4 +\|y\|^2 )  +C$ with $c_5 <2$, then it admits a $2\pi$-periodic probability solution.
\end{example}

Indeed, taking $a:=1$ and $ D:=2$, we have
\begin{eqnarray*}
\partial_t V(x,t) &=& \cos t\, (1- \exp{(-\|x\|^2) } ), \\
\nabla_x V(x,t) &=&  2(2+ \sin t)\,\exp{(-\|x\|^2) } \, x ,  \\
\nabla F(x)- x &=& \frac12 \exp{(\|x\|^2) }\, ( \|x\|^2 +1)\, x -x.
\end{eqnarray*}
In addition, it follows that $\|\nabla_x V(x,t) \| \leq 3$ and
$$
D^2 F(x)=\frac12 \exp(\|x\|^2) \,(\|x\|^2 +1) I_n + (\|x\|^2 +2 ) \exp{(\|x\|^2) } \, xx^{\top}  $$
is positive definite.
Moreover, it yields
\begin{eqnarray}
&&-\< \nabla_x V(x,t), \nabla F(x)- x \>  \nonumber\\
\label{e:Exam3}
&=& -(2+ \sin t)\, \big[ \|x\|^4 + \|x\|^2 -2 \|x\|^2\, \exp{(-\|x\|^2) }\big]  \nonumber\\
&\leq& - \|x\|^4 +3 .
\end{eqnarray}
Obviously, Hypotheses {\bf(H1)}-{\bf(H5)} hold with $ b=1, \, m=2, $ and $ M=3$.
Then  it yields from Theorem \ref{UF} that there exists a $2\pi$-periodic probability solution.
\qed

\begin{remark}
In Examples {\rm \ref{ex:4.2}-\ref{ex:4.4}}, the same conclusions hold by Theorem {\rm\ref{UF}} if the friction functions are replaced by $ F(x)= \| x\|^4 -|x\|^2 $ and $F(x)= \frac14 \int_0^{ \|x\|^2} e^s\, (s-1)\, ds$, respectively.
Meanwhile, these correspond to stochastic nonlinear relaxation oscillation systems.
\end{remark}
\bigskip
\section {Conclusions and Further Perspectives}\label{S:5}

This paper proves the  existence of a distributed periodic solution for stochastic time-periodic Newton systems under the assumptions that both the friction and potential functions tend to positive infinity at infinity, and the friction-potential gradient inner product grows at least like an even-power polynomial. This result is robust regarding bounded perturbations of the Hessian-friction matrix, at the cost of a slight increase in the growth rate of the potential function at infinity. Consequently, our method can be applied to systems like the one in Eq. \eqref{e:2.11}, where $V(x)=\|x\|^{6}$ and $D^2F(x)$ is replaced by $D^2\left(\|x\|^{4}\pm \|x\|^{2}\right)+B(x,y,t)$ with $B$ being bounded. To the best of our knowledge, none of the existing results can be applied to this system. Our results not only largely confirm the stochastic Levinson conjecture proposed in \cite{Liyong} under broad conditions and resolve the open problem and the interesting problem posed in \cite [P. 342]{Liyong}, but also possess advantages 1-6 enumerated in the Introduction, which are absent in existing results.

Our further research plan involves studying the uniqueness, exponential ergodicity, and exponential mixing of distributed time-recurrent solutions to the stochastic time-recurrent Newtonian equations, including distributed periodic solutions and quasi-periodic solutions, etc.

\bigskip
\noindent{\bf Acknowledgements} \quad The project is supported by the National Natural Science Foundation of China with
grant numbers 12171321.
\medskip

\noindent {\bf Data Availability} \qquad Our manuscript has no associate data.\\
\medskip
\\
{\bf \Large Declarations}
\smallskip\\
\noindent {\bf Conflicts of Interest}  \,\, All authors have no Conflict of interest.

\end{document}